\documentclass[11pt]{article} 
\usepackage{amssymb}
\usepackage{latexsym,amsmath,amsthm,amsfonts}

\newtheorem{prop}{Proposition}
\newtheorem{lem}{Lemma}
\newtheorem{thm}{Theorem}
\newtheorem{cor}{Corollary}

\newtheorem{lemref}{Lemma}

\newcommand{\be}{\begin{equation}}
\newcommand{\ee}{\end{equation}}

\def\real{\mathbb R}

\def\P{\mathbb P}

\newcommand{\deq}{\stackrel{\scriptscriptstyle\triangle}{=}}

\newcommand{\C}{\mathcal{C}}
\newcommand{\sigS}{\mbox{$\cal S$}}
\newcommand{\sigB}{\mbox{$\cal B$}}
\newcommand{\D}{\mbox{$\cal D$}}

\newcommand{\F}{\mbox{$\mathcal F$}}
\newcommand{\G}{\mbox{$\cal G$}}
\newcommand{\X}{\mbox{$\cal X$}}

\newcommand{\bX}{\mbox{$\bf X$}}
\newcommand{\bY}{\mbox{$\bf Y$}}

\begin{document}

\bibliographystyle{plain}

\title{The Gap Dimension and Uniform Laws of Large Numbers for
Ergodic Processes}

\author{Terrence M.\ Adams
\thanks{Terrence Adams is with the Department of Defense, 
9800 Savage Rd. Suite 6513, Ft. Meade, MD 20755} \ 
and Andrew B.\ Nobel
\thanks{Andrew Nobel is with the Department of Statistics and 
Operations Research, University of North Carolina, Chapel Hill,
NC 27599-3260.  Corresponding author.  Email: nobel@email.unc.edu}}

\date{May 24, 2010}

\maketitle

\begin{abstract}
Let $\F$ be a family of Borel measurable functions
on a complete separable metric space.
The gap (or fat-shattering) dimension of $\F$ is a combinatorial
quantity that measures the extent to which functions $f \in \F$ 
can separate finite sets of points at a 
predefined resolution $\gamma > 0$. 
We establish a connection between the gap dimension
of $\F$ and the uniform convergence of its sample
averages under ergodic sampling.  In particular, we show that if the
gap dimension of $\F$ at resolution $\gamma > 0$ is finite, 
then for every ergodic process the sample averages of functions
in $\F$ are eventually within $10 \gamma$ of their limiting 
expectations {\em uniformly} over the class $\F$.  If the gap dimension
of $\F$ is finite for every resolution $\gamma > 0$ then 
the sample averages of functions in $\F$ converge uniformly 
to their limiting expectations.  We assume only that $\F$ is uniformly
bounded and countable (or countably approximable).   No smoothness
conditions are placed on $\F$, and no
assumptions beyond ergodicity are placed on the 
sampling processes.  Our results extend existing work
for i.i.d.\ processes.

\end{abstract}

\vskip.5in

\noindent
Running Title:  Gap Dimension and Uniform Laws of Large Numbers

\vskip.5in


\newpage

\section{Introduction}

Let $\X$ be a complete separable metric space, and let
$\F$ be a countable family of Borel-measurable functions
$f: \X \to \real$.  We assume in what follows that 
$\F$ is uniformly bounded in the sense that 
$|f(x)| \leq M$ for every $x \in \X$ 
and $f \in \F$, where $M < \infty$ is a fixed constant.  
Let $\bX = X_1, X_2, \ldots$ 
be a stationary ergodic process taking values in $\X$.
By the ergodic theorem, for each $f \in \F$,
the averages
$m^{-1} \sum_{i=1}^m f(X_i)$ converges with probability 
one to $Ef(X)$.  Of interest here is the 
limiting behavior of the discrepancy
\be
\label{unif-discrep}
\Gamma_m(\F:\bX)
\ \deq \ 
\sup_{f \in {\cal F}} 
\left| \frac{1}{m} \sum_{i=1}^{m} f(X_i) \, - \, Ef(X) \, \right| ,
\ee
which measures the maximum difference between $m$-sample averages
and their limiting expectations over the functions in $\F$.

The discrepancy $\Gamma_m(\F:\bX)$ and 
related quantities have been studied 
in a number of fields, including empirical process theory,
machine learning and non-parametric inference.  
The majority of existing work considers the case in which 
$X_1, X_2, \ldots$ are independent and identically distributed, but
there is also a substantial literature concerned with
the behavior of the discrepancy for mixing
processes (see \cite{AdNob09} and the discussion below).
Our focus here is on the general dependent case: the process
$\bX$ is not assumed to satisfy any mixing conditions beyond
ergodicity.

When $\bX$ is ergodic, the limiting behavior of the discrepancy
$\Gamma_m(\F:\bX)$ can be summarized by a single
number.  As shown in Steele
\cite{Stee78}, Kingman's subadditive ergodic 
theorem implies that there is a non-negative constant 
$\Gamma(\F:\bX)$ such that
\be
\label{Steele}
\lim_{m \to \infty}
\Gamma_m(\F:\bX) \to
\Gamma(\F:\bX) \ \ \mbox{wp1} .
\ee
We will call $\Gamma(\F:\bX)$ the {\em asymptotic discrepancy} of $\F$
on $\bX$, and will omit mention of $\bX$ when no confusion
will arise.  When $\Gamma(\F:\bX) = 0$ the sample averages 
of function $f \in \F$ converge uniformly to their limiting expectations,
and $\F$ is said to be a Glivenko Cantelli class for 
the process $\bX$.  

In this paper we provide bounds on the asymptotic discrepancy of $\F$ in terms 
of a combinatorial quantity known as the gap dimension 
that measures the complexity of $\F$ at 
different resolutions or scales.  

\vskip.1in

\noindent
{\bf Definition:} Let $\gamma > 0$.  The family $\F$ is said to 
$\gamma$-shatter a finite set $D \subseteq \X$ if there
is an $\alpha \in \real$ such that for every $D_0 \subseteq D$ there exists a
function $f \in \F$ satisfying
\[
f(x) > \alpha + \gamma \mbox{ if } x \in D_0
\ \ \ \mbox{and} \ \ \ 
f(x) < \alpha - \gamma \mbox{ if } x \in D \setminus D_0
\]
The gap dimension of $\F$ at resolution $\gamma$, written $\dim_\gamma(\F)$, 
is the largest $k$ such
that $\F$ $\gamma$-shatters {\em some} set of cardinality $k$.  If $\F$ 
can $\gamma$-shatter sets of arbitrarily large finite cardinality, then
$\dim_\gamma(\F) = + \infty$.

\vskip.1in

The gap dimension was introduced by Kearns and Schapire \cite{KeaSch94}
in a slightly more general form.  Specifically, they allowed the constant $\gamma$
to be replaced by a fixed function $g: \X \to \real$.  We will refer to this
notion as the weak gap dimension in what follows.
The definition of gap dimension given here was suggested
by Alon, Ben-David, Cesa-Bianchi and Haussler \cite{AlBDCBHa97},
who also established elementary bounds relating the gap and
weak gap dimensions.  
Gap dimensions have been referred to by a variety of 
names in the literature,
including scale-sensitive dimension and fat-shattering dimension.
Our principal result is the following theorem.  As above, $\X$ is
assumed to be a complete separable metric space.

\begin{thm}
\label{FS1}
Let $\F$ be a countable, uniformly bounded family of Borel measurable 
functions $f: \X \to \real$, and let $\bX$ be a stationary ergodic process 
with values in $\X$.  If the asymptotic discrepancy 
$\Gamma(\F : \bX) > \eta$ for some $\eta > 0$,
then $\dim_\gamma(\F) = \infty$ for every $\gamma \leq \eta / 10$.
\end{thm}

The constant $10$ dividing $\gamma$ can, with
minor modifications of the proof, be improved to $4 + \epsilon$,
where $\epsilon$ is any fixed positive constant.
Theorem \ref{FS1} has the following, equivalent, form.

\begin{cor}
\label{CP}
Let $\F$ be as in Theorem \ref{FS1}.  If 
$\dim_\gamma(\F) < \infty$ for some $\gamma > 0$
then $\Gamma(\F : \bX) \, \leq \, 10 \gamma$
for every stationary ergodic process.
In particular, if $\dim_\gamma(\F) < \infty$ 
for every $\gamma > 0$, then
$\Gamma(\F : \bX) = 0$ for every stationary ergodic
process.
\end{cor}

\paragraph{Uncountable Families} 

The countability of $\F$ ensures that the
discrepancies $\Gamma_m(\F, \bX)$, $m \geq 1$, are measurable. 
More importantly, countability of $\F$ is  
used in the proof of Proposition \ref{Tree1} and 
is a key assumption in Lemma \ref{RedLem3}. 
Nevertheless, one may readily extend Theorem \ref{FS1}
to uncountable families under simple approximation
conditions.  Call a (possibly uncountable) family $\F$ nice for a process
$\bX$ if $\Gamma_m(\F:\bX)$ is measurable for each $m \geq 1$, and if
for every $\epsilon > 0$ there exists a countable sub-family
$\F_0 \subseteq \F$ such that 
$\limsup_m \Gamma_m(\F:\bX) \leq 
\limsup_m \Gamma_m(\F_0 : \bX) + \epsilon$ with probability
one.  The conclusion of Theorem \ref{FS1} immediately 
extends to any ergodic processes $\bX$ for which $\F$ is nice. 


In spite of such extensions, assumptions regarding the countability 
or countable approximability of $\F$ cannot be dropped
altogether, as they exclude extreme examples that
can arise in the context of dependent processes.  
We illustrate with a simple example from
\cite{AdNob09}.  Let $T$ be an irrational rotation of the unit circle
$S_1$ with its uniform measure.  Denote by $T^i$ the $i$-fold 
composition of $T$ with itself if $i \geq 1$, the  $i$-fold 
composition of $T^{-1}$ with itself if $i \leq -1$ and the
identity if $i = 0$.
For each $x \in S_1$ let 
$\C_{x} = \cup_{i=-\infty}^{\infty} \{ T^i x \}$
be the (bi-infinite) trajectory of $x$ under $T$, and
let $\F$ be the family of indicator functions of the 
sets $\C_x$.  Note that $\F$ is uncountable, and
that every set $\C_x$ has Lebesgue measure zero.
For distinct points $x_1, x_2 \in S_1$, either 
$C_{x_1} = C_{x_2}$, or 
$C_{x_1} \cap C_{x_2} = \emptyset$, and therefore
$\dim_\gamma(\F) = 1$ for $0 < \gamma < \frac{1}{2}$. 
Now let $X_i = T^i X_0$, where $X_0$ is uniformly distributed
on $S_1$.  Then the process $\bX = X_1, X_2, \ldots$
is stationary and ergodic.  Moreover, it is easy to see that $Ef(X) = 0$
for each $f \in \F$, and that 
$\sup_{f \in {\cal F}}  m^{-1} \sum_{i=1}^m f(X_i) = 1$.
Thus $\Gamma_m(\F : \bX) = 1$ with probability one for
each $m \geq 1$, and the conclusion of Corollary \ref{CP}
fails to hold.

\subsection{Related Work}

Vapnik and Chervonenkis \cite{VC81} gave necessary and sufficient conditions
for uniform convergence of sample means in the i.i.d.\ case. 
Specifically, they showed that if $\bX$ is i.i.d., then 
$\Gamma(\F: \bX) = 0$ if and only if
$n^{-1} \log N(\epsilon, \F, X_1^n) \to 0$ in probability
for every $\epsilon > 0$.  Here $N(\epsilon, \F, X_1^n)$
is the number of $\epsilon$-balls needed to cover $\F$ under the
empirical $L_1$ metric 
$d(f_1,f_2) = n^{-1} \sum_{i=1}^n |f_1(X_i) - f_2(X_i)|$.
Extensions of these results to empirical processes can be found,
for example, in Gin\'e and Zinn \cite{GinZin84} (see also Dudley
\cite{Dudley99}).

Talagrand \cite{Tala87} gave necessary and sufficient conditions for 
uniform convergence of sample means, which are different than those of
\cite{VC81}. 
He showed that $\Gamma(\F : \bX)  > 0$
for an i.i.d.\ process $\bX$ with $X_i \sim P$ if and only if there exists a set
$A$ with $P(A) > 0$ and $\gamma> 0$ such that for every $n \geq 1$
the family $\F$ $\gamma$-shatters $P^n$-almost every sequence 
$x_1,\ldots, x_n \in A^n$.

Alon {\em et al.} \cite{AlBDCBHa97} considered the relationship between the gap
dimension and the learnability of classes of uniformly bounded
functions under independent sampling.  In particular, they
showed that if $\F$ is a family of functions $f : \X \to [0,1]$ 
satisfying suitable measurability conditions, and such that
$\dim_\gamma(\F)$ is finite
for some $\gamma >0$, then
\be
\label{ACBH}
\lim_{n \to \infty} 
\left[
\sup_{{\bf X} \in {\cal I}({\cal X})}
\P \left( \sup_{m \geq n} \Gamma_m(\F : \bX) > \varepsilon  \right)
\right]
\ = \ 0
\ee
when $\varepsilon = 48 \gamma$.  
Here ${\cal I}({\cal X})$ is the family of all i.i.d.\ 
processes taking values in $\X$.  Conversely, if 
$\dim_\gamma(\F) = +\infty$, they showed that (\ref{ACBH}) fails
to hold for every $\varepsilon < 2 \gamma$.  Further
connections between the gap dimension and different notions
of learnability (in the i.i.d.\ case) can be found in
\cite{BaLoWi94} and the references therein.
Talagrand \cite{Tala03} and Mendelson and Vershynin \cite{MenVer03}
showed that the $L_2$ covering numbers of a uniformly
bounded sets of functions can be bounded in terms of its
weak gap dimension.

In addition to the papers cited above, there are a number of results 
on uniform convergence for dependent 
processes satisfying a variety of standard mixing conditions;
a discussion of these results can be found in \cite{AdNob09}.
In related work,
Rao \cite{Rao62} and Billingsley and Tops\o e \cite{BilTop67} studied and
characterized classes of functions $\F$ such that
$\sup_{\cal F} | \int f dP_n  - \int dP| \to 0$ whenever $P_n$
converges weakly to $P$.  As noted in \cite{BilTop67}, the elements
of such uniformity classes are necessarily continuous almost everywhere with
respect to $P$.  Bickel and Millar \cite{BicMil92} provided sufficient
conditions for a more general notion of uniformity, and revisited several
of the results in earlier papers.

Adams and Nobel \cite{AdNob09} established Theorem \ref{FS1} in
the special case where the elements of $\F$ are indicator
functions of subsets of $\X$.  The problem
simplifies in this case, as $\dim_\gamma(\F)$ is zero
for $\gamma \geq 1/2$, and equal to the VC-dimension
of $\F$ if $0 \leq \gamma < 1/2$.  If
$\F$ has finite VC-dimension, their results imply that
$\Gamma(\F: \bX) = 0$ for every ergodic process $\bX$. 
For uniformly bounded families $\F$ they show that
$\Gamma(\F: \bX) = 0$ for every ergodic process $\bX$ if 
$\dim_0(\F) < \infty$, or if $\F$ is a VC-graph class (c.f.\ \cite{Poll84}).

\subsection{Overview}

The proof of Theorem \ref{FS1} is based on the direct 
construction of $\gamma$-shattered sets of arbitrarily large cardinality.
In particular, the proof does note make use of results or techniques 
from the study of uniform convergence in the i.i.d.\ case. 
The core of the construction, which is contained in 
Section \ref{PfTree1} below, 
follows the arguments in \cite{AdNob09}.  

In the next section we reduce Theorem \ref{FS1}
to an analogous result with $\X$ is equal to the unit 
interval.  This equivalent result is stated in Theorem \ref{FS2}.
Section \ref{Pre} contains several preliminary definitions and
Lemmas used in the proof of Theorem \ref{FS2}.
The proof of Theorem \ref{FS2} is presented in Sections \ref{OutT2} - 
\ref{PfFS2}.  Section \ref{OutT2} gives an outline of the proof of
the theorem.  The proofs of two key propositions are given in
Sections \ref{PfTree1} and \ref{PfFullJoin}.
The diagram below provides an overview of the proof.
\begin{eqnarray*}
\mbox{Theorem \ref{FS1}} 
\ \ \Leftarrow \ \  
\mbox{Theorem \ref{FS2}}
& \Leftarrow &  
\mbox{Proposition \ref{FullJoin}} \ + \ \mbox{Lemma \ref{Join}}
\ + \ \mbox{Lemma \ref{RedLem3}} \\
&  & 
\hspace*{.39in} \Uparrow \\
& &
\mbox{Proposition \ref{Tree1}} \ + \ \mbox{Lemma \ref{Ptree}} \\
\end{eqnarray*}

\section{Reduction to the Unit Interval}
\label{RUI}

Let $\X$ and $\F$ be as in Theorem \ref{FS1} and let $\bX$ be an 
$\X$-valued ergodic
process, defined on an underlying probability space
$(\Omega, {\mathcal A}, \P)$,
such that $\Gamma(\F, \bX) > \eta > 0$.
By assumption, there exists a number $0 < M < \infty$ such 
that $|f| \leq M$ for each $f \in \F$.   Replacing $f \in \F$ with
$f' = (f + M) / 2M$, we may 
assume without loss of generality that each $f \in \F$ takes values 
in $[0,1]$.  The proof
of the following lemma, which relies on elementary ergodic
theory, is similar to that of Lemma 5 in \cite{AdNob09}, and is omitted.

\begin{lemref}
\label{RedLem1}
Let $\bX$ be a stationary ergodic process with values in $\X$. 
If $\Gamma(\F : \bX) > \eta > 0$, then $\X$ is necessarily uncountable, 
and there exists a stationary ergodic process
$\tilde{\bX}$ with values in $\X$
such that $\P(\tilde{X}_i = x) = 0$ for each $x \in \X$
and $\Gamma(\F: \tilde{\bX}) > \eta$.
\end{lemref}

Let $\mu(\cdot)$ be the marginal distribution of $\bX$. 
By Lemma \ref{RedLem1}, it suffices to establish Theorem \ref{FS1} in
the case where $\X$ is uncountable, and $\mu(\cdot)$ is non-atomic.
Let $\lambda(\cdot)$ denote ordinary Lebesgue measure on the unit
interval $[0,1]$ equipped with its Borel subsets $\sigB$.
By standard results in real analysis ({\em c.f.} Theorem 5.16 of \cite{Royd88}), 
there is a measure space
isomorphism between $(\X, \sigS, \mu)$ and $([0,1], \sigB, \lambda)$.
More precisely, there exist Borel measurable 
sets $\X_0 \subseteq \X$ and $I_0 \subseteq [0,1]$,
and a bijection $\psi: \X_0 \to I_0$ with the following properties:
(i) $\mu(\X_0) = \lambda(I_0) = 1$; (ii)
$\psi$ and $\psi^{-1}$ are measurable with respect to the restricted sigma
algebras $\sigS \cap \X_0$ and $\sigB \cap I_0$, respectively; and (iii)
$\mu(A) = \lambda(\psi(A))$ for each $A \in \sigS \cap \X_0$.
In particular, the event $E = \{ X_i \in \X_0^c \mbox{ for some } i \geq 1 \}$ 
has probability zero.  By removing $E$ from the underlying sample space, 
we may assume without loss of generality
that $X_i(\omega) \in \X_0$ for each sample point 
$\omega$ and each $i \geq 1$.  

Define $Y_i = \psi(X_i)$ for $i \geq 1$.  Then the process
$\bY = Y_1, Y_2, \ldots \in [0,1]$ is stationary and ergodic with 
marginal distribution $\lambda$.  For each function $f \in \F$ define
an associated function $\tilde{f} : [0,1] \to [0,1]$ via the rule
\[  
\tilde{f} \ = \ 
\left\{ \begin{array}{ll} 
          (f \circ \psi^{-1}) (u) & \mbox{if $u \in I_0$} \\ 
          0 & \mbox{otherwise} 
               \end{array} \right. 
\]
and let $\tilde{\F} = \{ \tilde{f} : f \in \F \}$.  
It is easy to see that $\tilde{f}(Y_i) = f(X_i)$, and in particular,
that $E \tilde{f}(Y) = Ef(X)$.  Thus
$\Gamma_m(\tilde{\F} : \bY) = \Gamma_m(\F : \bX)$ with probability
one for each $m \geq 1$.
Moreover, if $k$ distinct points 
$u_1 ,\ldots, u_k \in [0,1]$ are $\gamma$-shattered by $\tilde{\F}$, 
then necessarily each $u_j \in I_0$, and
the (distinct) points $\psi^{-1}(u_1), \ldots, \psi^{-1}(u_k) \in \X$
are $\gamma$-shattered by $\F$.
It follows that $\dim_\gamma(\tilde{\F}) \leq \dim_\gamma(\F)$.
Theorem \ref{FS1} is therefore a corollary of the following result.

\begin{thm}
\label{FS2}
Let $\F$ be a countable family of Borel measurable 
functions $f: [0,1] \to [0,1]$, and let $\bX = X_1, X_2, \ldots \in [0,1]$ be 
a stationary ergodic process with $X_i \sim \lambda$.  
If the asymptotic discrepency $\Gamma (\F : \bX) > \eta > 0$
then $\dim_\gamma(\F) = \infty$ for every $\gamma \leq \eta / 10$.
\end{thm}

\section{Preliminaries}
\label{Pre}

In this section we define three elementary notions that will be
used in the proof of Theorem \ref{FS2}.  The first is the segments of a 
function $f: [0,1] \to [0,1]$.  The second is the join of a sequence of 
families of disjoint sets.  The third is an ancestral set in a binary
tree.  Lemma \ref{Join} establishes a simple connection between
joins, segments and the gap dimension.  Lemma \ref{Ptree} provides
a useful bound for obtaining a subtree with good ancestral
properties from a large initial binary tree.

\subsection{Segments and Regular Families}
\label{SRF}

Let $\F$ and $\bX$ be as in the statement of Theorem \ref{FS2},
and suppose that $\Gamma (\F : \bX) > \eta > 0$.  
Assume without loss of generality that $\eta$ is rational, and let
$\gamma = \eta / 5$. Let $K = \lfloor \gamma^{-1} \rfloor + 1$ if $\gamma^{-1}$ 
is not an integer, and $K = \gamma^{-1}$ otherwise.
For each $f \in \F$ and $1 \leq k \leq K$ define sets
\be
\label{segdef}
s_k(f) \ = \ 
\left\{ \begin{array}{ll} 
          f^{-1} [(k-1) \, \gamma, k \, \gamma) & \mbox{if \ $1 \leq k \leq K-1$} \\ [.06in]
          f^{-1} [(K-1) \, \gamma, 1] & \mbox{if \ $k = K$.} 
               \end{array} \right. 
\ee

\vskip.15in

\noindent
{\bf Definition:}
The sets $s_k(f)$ will be called {\em $\gamma$-segments} of $f$.  
Let $\pi(f) = \{ s_k(f) : 1 \leq k \leq K \}$ be the partition of $[0,1]$
generated by the $\gamma$-segments of $f$.
Two segments
$s_k(f)$ and $s_{k'}(f)$ will be called adjacent if they correspond to adjacent 
intervals, equivalently if $|k - k'| =1$, and non-adjacent 
if $|k - k'| \geq 2$.  

\vskip.15in

In order to establish Theorem \ref{FS2}, we first consider 
families $\F$ whose elements
satisfy a topological regularity condition.  
Given a family $\F$ of functions $f:[0,1] \to [0,1]$, define the
associated collection of sets
\be
\label{regsets}
\C(\F) \ = \ \{ f^{-1} [a,b) :  0 \leq a < b < 2 \mbox{ rational, and }  f \in \F \} .
\ee
Including values $b > 1$ ensures that $\C(\F)$ contains sets of the form
$f^{-1} [a,1]$.  Note that $\C(\F)$ is countable if $\F$ is countable.  

\vskip.13in

\noindent
{\bf Definition:}
A family $\F$ of measurable functions $f : [0,1] \to [0,1]$
is {\em regular} if it is countable, and each element of 
$\C(\F)$ is a finite union of intervals.

\subsection{Joins and the Gap Dimension}

In ergodic theory, the join of a finite collection of sets contains the atoms of 
their generated field.  Here we employ a minor generalization of this 
notion.

\vskip.1in

\noindent
{\bf Definition:} Let $\D_1, \ldots, \D_k$ be finite families 
sets in $[0,1]$ such that the elements of each family are disjoint.
The join of $\D_1,\ldots, \D_k$,
denoted $\bigvee_{i=1}^k \D_i$ or $\D_1 \vee \cdots \vee \D_k$,
is the collection of all 
{\em non empty} intersections $D_1 \cap \cdots \cap D_k$ 
where $D_i \in \D_i$ for $i = 1,\ldots,k$. 

\vskip.14in

The next lemma establishes a useful connection between 
the gap dimension of $\F$ and the join of non-adjacent
segments of functions $f \in \F$.  Its proof is based 
on similar results in \cite{Mat02} and \cite{AdNob09}.

\begin{lem}
\label{Join}
Suppose that for some $L \geq 1$
there exists a sub-family $\F_0 \subseteq \F$ of $2^L$
functions, and a pair $k,k' \in [K]$
of non-adjacent integers such that the join
\[
J \ = \ \bigvee_{f \in {\cal F}_0} \{ s_k(f), s_{k'}(f) \} 
\]
of non-adjacent $\gamma$-segments
has cardinality $2^{2^L}$.  Then $\dim_{\gamma/2}(\F) \geq L$.
\end{lem}

\noindent
{\bf Remark:}
The conditions of the lemma ensure that each of the possible
intersections contained in $J$ is non-empty, and therefore $J$ has
maximum cardinality.  

\vskip.1in

\noindent
{\bf Proof:} Indexing the elements of $\F_0$ in an arbitrary manner by
subsets of $[L] : = \{1,\ldots,L\}$, we may write 
$\F_0 = \{ f_\alpha : \alpha \subseteq [L] \}$.  For $i = 1,\ldots,L$, let $x_i$ be
any element of the intersection
\[
\left( \bigcap_{ \alpha \subseteq [L], i \in \alpha } s_k(f_\alpha) \right)
\, \cap \,
\left( \bigcap_{ \alpha \subseteq [L], i \not\in \alpha } s_{k'}(f_\alpha) \right) ,
\]
which is non-empty by assumption.  Suppose without loss of generality
that $k < k'$, and let $c = \gamma (k + k' - 1)/ 2$.  Let $\beta$ be any subset of $[L]$
and consider the corresponding function $f_\beta \in \F_0$.  
If $i \in \beta$, the selection of $x_i$ ensures that 
$x_i \in s_k(f_\beta)$, and consequently 
$f_\beta(x_i) < \gamma k < c - \gamma / 2$.  On the other
hand, if $i \in \beta^c$ then $x_i \in s_{k'}(f_\beta)$, and in this case 
$f_\beta(x_i) \geq \gamma (k' - 1) \geq c + \gamma / 2$.
As $\beta$ was arbitrary, it follows that $\dim_{\gamma/2}(\F) \geq L$.

\subsection{Binary Trees and Ancestral Sets}

Binary trees appear in several key results of the paper.  
Throughout we consider standard binary trees $T$ that 
have a single root, which is assumed to be  
located at the top of the tree.  Vertices of $T$ are referred to
as nodes, and usually denoted by $s$ or $t$.  Each node of $T$ has either zero or 
two distinct children and, with the exception of the root, a single parent.  A node with
two children is said to be internal; a node with no children is called a leaf.
The set of leaves in a tree $T$ will be denoted by $\tilde{T}$.
A descending path in $T$ is a sequence of adjacent nodes that 
proceeds only from parent to child.  The depth, or level, of a node $t \in T$ 
is the length
of the shortest (necessarily descending) path from the root to $t$.  
The set of nodes at level $r$ of $T$ will be denoted $T[r]$.  The depth of 
$T$ is the maximum depth of any node in $T$.  We will exclusively consider 
trees of finite depth, say $L$, that are complete in the sense that
$T[r]$ contains $2^r$ nodes for $r = 0, \ldots, L$.  In this case,
$\tilde{T} = T[L]$ and each node $t \in T[r]$ with $0 \leq r \leq L-1$ is
internal.

\vskip.1in

\noindent
{\bf Definition:}
Let $T$ be a binary tree.  
A node $s$ in $T$ is an ancestor of a node $t$ if there is a 
descending path in $T$ from $s$ to $t$ of length greater than
or equal to one.
A node $s$ will be called an ancestor 
of a set $A \subseteq T$ if $s$ is an ancestor of 
{\em some} $t \in A$.

\vskip.1in

The next Lemma establishes a pigeon-hole type 
result showing that
any large collection of leaves must have a 
correspondingly large set of ancestors
in some nearby level of the tree.

\begin{lem}
\label{Ptree}
Let $T$ be a full binary tree of depth $L$, and let $\tilde{T}$ denote
the $2^L$ leaves of $T$.  Suppose that there exists a set of leaves
$S \subset \tilde{T}$ and a constant $0 < c < 1$ such that $|S| \geq c 2^L \geq 4$. 
Let $u = \lceil \log_2 c^{-1} + 1 \rceil$. 
Then there exists a set $S' \subseteq T[l_0]$ with $L-u \leq l_0 \leq L-1$ 
such that for each node $s \in S'$ both of its children are ancestors 
of $S$, and 
\be
\label{blbd}
|S'| 
\ \geq \ 
\frac{c 2^L}{4L} .
\ee
\end{lem}

\noindent
{\bf Proof:} For $l = 1,\ldots, L-1$, let $m_l$ be the number of nodes $s$ at level
$l$ that are the ancestor of some node $t \in S$, and let $n_l$ be the
number of nodes at level $l$ with the property that both their children are ancestors
of a node $t \in S$.  It is easy to see that $|S| = m_{L-1} + n_{L-1}$, and more 
generally we have
\[
|S| 
\ = \ 
m_{L-v} \, + \, n_{L-v} \, + \, n_{L-v+1} \, + \cdots + \, n_{L-1}
\ \leq \ 
2^{L-v} \, + \, \sum_{l=L-v}^{L-1} n_l
\]
for $v = 1,\ldots, L-1$.  Setting $v = u$, the
assumption that $|S| \geq c 2^L$ yields
\[
\sum_{l=L-u}^{L-1} n_l
\ \geq \ 
c 2^L \, - \, 2^{L-u}
\ = \ 
2^{L - u} ( c 2^u - 1)
\ \geq \ 
2^{L - u} ,
\]
where the last inequality follows from the definition of $u$.  
Let $n_{l_0}$ be the largest value of $n_l$ appearing in the sum
above, and let $S'$
be the nodes at level $l_0$ of $T$ with the property that both 
their children are ancestors of $S$.  Then 
\[
|S'| 
\ = \ 
n_{l_0}
\ \geq \ 
\frac{2^{L - u}}{u} 
\ \geq \ 
\frac{c 2^L}{4u} 
\ \geq \ 
\frac{c 2^L}{4L} 
\]
where the second inequality follows from the definition of $u$.

\section{Outline of the Proof of Theorem \ref{FS2}}
\label{OutT2}

In this section we present an outline of the proof of 
Theorem \ref{FS2}.  We begin with Proposition \ref{Tree1},
which is the key result of the paper.  The proposition
shows that if $\F$ is regular and $\Gamma(\F : \bX) > 0$ then
one can associate the nodes of an arbitrarily large binary tree with
segments of select functions in $\F$ 
in such a way that 
(i) the intersection of segments along every path from the 
root to a leaf is non-empty, and (ii) sibling segments are 
non-adjacent.  The resulting structure will be called
an intersection tree.

Proposition \ref{FullJoin} refines Proposition \ref{Tree1} 
using the pigeon-hole principle from Lemma \ref{Ptree}.  
It ensures that for every finite $L \geq 1$ there is
a family of $L$ functions in $\F$ having non-adjacent
segments with maximal join.
The final step in the proof of Theorem \ref{FS2} is to
remove the regularity condition on $\F$.  This is done by
means of a measure space isomorphism described in 
Lemma \ref{RedLem3}.  The proof of Theorem \ref{FS2}
appears in Section \ref{PfFS2}.

\subsection{Intersection Trees}

\begin{prop}
\label{Tree1}
Let $\F$ and $\bX$ be as in Theorem \ref{FS2}.
Suppose that $\Gamma (\F : \bX) > \eta > 0$ and that $\F$
is regular.
Then for each $L \geq 1$ there exists functions 
$g_1, \ldots, g_L \in \F$ and 
a complete binary tree $T$ of depth $L$
such that each node $t \in T$ is associated with a subset $B_t$ of $[0,1]$
in such a way that the following two conditions are satisfied.
\begin{enumerate}

\item[(a)]
For each internal node $t \in T$ at level $\ell$, the
sets $B_{t'}$ and $B_{t''}$ associated with its children
$t'$ and $t''$ are 
equal to non-adjacent segments of $g_{\ell+1}$.

\item[(b)]
For each node $t \in T$, the intersection $W_t$ of the sets $B_s$ 
appearing along a descending path from the root to $t$ 
has non-empty interior.

\end{enumerate}
\end{prop}

The proof of Proposition \ref{Tree1} is given in 
Section \ref{PfTree1}.

\subsection{Maximal Joins}

\begin{prop}
\label{FullJoin}
Let $\F$ and $\bX$ be as in Theorem \ref{FS2}.
Suppose that $\Gamma (\F : \bX) > \eta > 0$ and that $\F$ is
regular.  Let $\gamma = \eta/5$.
For each $L \geq 1$ there are functions
$f_1, \ldots, f_L \in \F$ and a pair $k,k' \in [K]$ of non-adjacent integers
such that the join 
\[
J \ = \ \{ s_k(f_1), s_{k'}(f_1) \} \vee \cdots \vee \{ s_k(f_{L}), s_{k'}(f_{L}) \}
\]
of non-adjacent $\gamma$-segments
has (maximum) cardinality $2^L$, and every element of $J$
has positive Lebesgue measure.  
\end{prop}

The proof of Proposition \ref{FullJoin}
appears in Section \ref{PfFullJoin} below.

\subsection{Removing Regularity}

Together, Lemma \ref{Join} and Proposition \ref{FullJoin} 
establishes Theorem \ref{FS2} in 
the special case of regular families.
In order to remove the assumption of regularity, 
we require the following result, whose proof can 
be found in \cite{AdNob09}.

\vskip.1in

\begin{lemref}
\label{RedLem3}
Let $\C =\{ C_1, C_2,\ldots \}$ be a countable collection of Borel 
subsets of $[0,1]$ such that the maximum diameter of the elements of the 
join $J_n = \bigvee_{i=1}^{n} \{ C_i, C_i^c \}$ tends to zero as $n \to \infty$.  
Then there exists a Borel-measurable map $\phi: [0,1] \to [0,1]$ 
and a Borel set $V_1\subseteq [0,1]$ of measure one
such that: 
(i) $\phi$ preserves Lebesgue measure and is 1:1 on $V_1$;
(ii) the image $V_2 = \phi(V_1)$ and the inverse map 
$\phi^{-1} : V_2 \to V_1$ are Borel measurable;
(iii) $\phi^{-1}$ preserves Lebesgue measure; and
(iv) for every set $C \in \C$ there is a set $U(C)$, equal to a finite
union of intervals, such that 
$\lambda (\phi(C) \triangle U(C)) = 0$, where $\triangle$ is the 
usual symmetric difference.
\end{lemref}

\vskip.1in

\noindent
{\bf Remark:}
Lemma \ref{RedLem3}
is applied to the family of sets $\C = \C(\F)$.  The existence of the 
isomorphism $\phi$ requires that $\C$ be countable,
and this leads to the requirement that $\F$ be countable as well.

\vskip.1in

The proof of Theorem \ref{FS2} is given in Section \ref{PfFS2} below.

\section{Proof of Proposition \ref{Tree1} }
\label{PfTree1}

Construction of the intersection tree in Proposition \ref{Tree1} is based 
on a multi-stage procedure that is detailed below.  
At the first stage, we
produce a refining sequence $J_1, J_2, \ldots$ of joins in $[0,1]$
and simultaneously identify a sequence of functions $f_1, f_2, \ldots \in \F$.
The join $J_n$ is generated from selected non-adjacent 
segments of $f_1,\ldots,f_n$.
The function $f_{n+1}$ chosen at step $(n+1)$ is an element of 
$\F$ whose average differs from its
expectation by at least $\eta$ on a sample sufficiently large to 
ensure that the relative frequency of every element $A \in J_n$ is close
to its probability.
From $J_n$ and $f_{n+1}$ we identify a set $G_n$ equal 
to the union of the cells in $J_n$ on which
the average of $f_{n+1}$ is far from its expectation.   
The sets $G_n$ are used, in turn,
to produce a limiting ``splitting'' set $R_1$ via a weak convergence argument.  
This sequential process
is repeated in subsequent stages, with the important feature 
that the splitting sets $R_1,\ldots, R_{s-1}$
identified at stages $1,\ldots,s-1$ are used to 
generate the joins and the splitting set at stage $s$.   

The proof of Proposition \ref{Tree1} follows
the proof of Proposition 3 in \cite{AdNob09}. 
The earlier proposition treats 
the special case in which the
elements of $\F$ are indicator functions of sets, and hence binary
valued.  
The definition and construction of the splitting sets $R_s$ follow 
the arguments in the binary case, the principal difference being that the
generalized joins 
defined here involve segments rather than sets.  
The proof of Lemma \ref{Alem} below and the three
displays preceding it are identical
to arguments in \cite{AdNob09}.  Differences 
in the proofs emerge from the focus here on 
non-adjacent segments. 
In particular, the use of intersection trees 
or a similar hierarchical structure appears to be required, 
and the arguments that follow  
Lemma \ref{Alem} are somewhat more involved than in the binary case.  

The proof of Proposition \ref{Tree1} requires that one carefully
keep track of the quantities appearing at each step 
and stage of the construction, and how these quantities are defined.
For this reason, and due to the differences discussed above,
it is not possible to substantially shorten 
the proof Proposition \ref{Tree1} by an appeal to the earlier results.  
We provide a detailed argument below for completeness.

\subsection{Initial Construction}

Let $\F$ be a countable family of Borel measurable 
functions $f: [0,1] \to [0,1]$, and let $\bX = X_1, X_2, \ldots \in [0,1]$ be 
a stationary ergodic process defined on an underlying probability space
$(\Omega, {\cal A}, \P)$ such that $X_i \sim \lambda$. 
Assume that $\Gamma (\F : \bX) > \eta > 0$, and that every element
of $\C(\F)$ is a finite union of intervals.
Let $\delta = \eta/12$, and note that $0 < \delta < 1$.  For each $n \geq 1$ 
let
\[ 
\D_n = \{ [\, k \, 2^{-n}, (k+1) \, 2^{-n}) : 0 \leq k \leq 2^n-2 \} 
\cup \{ [1 - 2^{-n},1] \}
\]
be the $n$th order dyadic subintervals of $[0,1]$,
and let $\D = \cup_{n \geq 1} \D_n$.  
The set $A_0$ consisting of the endpoints of the intervals 
from which the elements of 
$\C(\F)$ and $\D$ are constructed is countable, and therefore 
has Lebesgue measure zero.  Removing a $\P$-null set of
outcomes from $\Omega$, we may assume that
$X_i(\omega) \in A_0^c$ for each $\omega \in \Omega$
and for every $i \geq 1$.  (This assumption is used
in the last part of the proof.)  

Below we identify a sequence of splitting sets 
$R_1, R_2, \ldots \subseteq [0,1]$ in stages, and
then use these sets to construct the intersection tree.

\vskip.1in

\noindent
{\bf Stage 1.}
The first stage of the construction proceeds as follows.
Let $f_1$ be any function in $\F$, and suppose 
that functions $f_1,\ldots,f_n \in \F$ have already been selected.
Let $J_n = \D_n \vee \pi(f_1) \vee \cdots \vee \pi(f_n)$ be
the join of the dyadic intervals of order $n$ and the $\gamma$-segments 
of the previously selected functions.  Here and in what follows
we take $\gamma = \eta/5$.
For each $\omega \in \Omega$, each function $g: [0,1] \to [0,1]$, 
and each $m \geq 1$, define the (pointwise) discrepancy
\be
\label{discrep}
\Delta^\omega(g : m)
\ \deq \ 
\left| \frac{1}{m} \sum_{i=1}^{m} g( X_i(\omega) ) - E g(X) \right| ,
\ee
which measures the difference between the expectation
of $g(X)$ and its average over the sample sequence 
$X_1(\omega),\ldots, X_m(\omega)$.
From the ergodic theorem and Proposition \ref{Steele}, 
it follows that there exists a sample point 
$\omega_{n+1} \in \Omega$, an integer 
$m_{n+1} \geq 1$ and a function $f_{n+1} \in \F$ such that 
\be
\label{inq1-init}
\Delta^{\omega_{n+1}} (I_A : m_{n+1}) \leq  \delta \, \lambda(A)
\mbox{ for each } A \in J_n 
\ee
and 
\be
\label{inq2-init}
\Delta^{\omega_{n+1}} (f_{n+1} : m_{n+1}) > \eta .
\ee
Defining the join
$J_{n+1} = \D_n \vee \pi(f_1) \vee \cdots \vee \pi(f_{n+1})$ 
and continuing, we may
select functions
$f_{n+2}, f_{n+3}, \ldots \in \F$ in a similar fashion.

The relations (\ref{inq1-init}) and (\ref{inq2-init}) 
together ensure that for many cells
$A \in J_n$ the average of $f_{n+1}$ on $A$ differs from its expectation
over $A$.  To make this precise, define the family
\[
H_{n}
\ = \ 
\left\{ A \in J_n \ : \  
   \Delta^{\omega_{n+1}}(f_{n+1} \cdot I_A : m_{n+1}) 
   \, > \, \frac{\eta}{2} \, \lambda(A) \right\} .
\]
As the next lemma shows, the sets in $H_{n} \subseteq J_n$ occupy a non-trivial
fraction of the unit interval.

\vskip.in

\begin{lem}
\label{Glem1}
If $G_{n} = \cup H_{n}$ is the union of the sets $A \in H_{n}$, then
$\lambda(G_{n}) \geq \eta / 6$.
\end{lem}

\noindent
{\bf Proof:} To simplify notation, let $\omega = \omega_{n+1}$, 
$f = f_{n+1}$, and $m = m_{n+1}$.
Decomposing $\Delta^\omega(f:m)$ over the elements
of $J_n$ and applying the triangle inequality, we obtain the bound
\[
\eta 
\ \leq \ 
\sum_{A \in H_{n}} \Delta^\omega(f \cdot I_A : m) 
\ + \ 
\sum_{A \in J_n \setminus H_{n}} \Delta^\omega(f \cdot I_A : m). 
\]
By definition of $H_{n}$, the second term is at most $\eta / 2$.  
The first term is at most
\begin{eqnarray*}
\lefteqn{
\sum_{A \in H_{n}} \Delta^\omega(f \cdot I_A : m) } \\ 
& \leq & 
\sum_{A \in H_{n}}
\left[ 
\frac{1}{m} \sum_{i=1}^{m} (f \cdot I_A) (X_i(\omega)) 
 \ + \,  E(f \cdot I_A)(X)
 \right] \\
& \leq & 
\sum_{A \in H_{n}} 
\left[
\frac{1}{m} \sum_{i=1}^{m} I_A(X_i(\omega))
\, + \,   
\lambda (A) 
\right] \\ [.1in]
& \leq &
\sum_{A \in H_{n}} 
\Delta^\omega(A : m)
\, + \,   
2 \, \lambda (G_{n}) \\ [.1in]
& \leq &
(\delta + 2) \lambda (G_{n})
\ \leq \ 
3 \, \lambda (G_{n}) .
\end{eqnarray*}
where the first inequality follows from the fact that $0 \leq f \leq 1$.
Combining the bounds above yields the stated inequality.

\vskip.1in

For each $n \geq 1$ define a sub-probability
measure $\lambda_n (B) = \lambda (B \cap G_n)$ on $([0,1],\sigB)$,
where $G_n = \cup H_{n}$.
The collection $\{ \lambda_n \}$ is tight, and is such that
$\lambda_n([0,1]) \geq \eta / 6$ for each $n$.  There is therefore a subsequence 
$n(1) < n(2) < \cdots$ such that $\lambda_{n(r)}$ converges weakly to a 
sub-probability measure $\nu_1$ on $([0,1],\sigB)$.
It is easy to see that $\nu_1$ is absolutely continuous with respect 
to $\lambda$, that $\nu_1([0,1]) \geq \eta / 6$, and that
the Radon-Nikodym derivative $d\nu_1 / d\lambda$ is 
is bounded above by 1.  
Define $R_1 = \{ x: (d\nu_1 / d\lambda)(x) > \delta \}$.  From the previous 
remarks it follows that 
\begin{eqnarray}
\frac{\eta}{6} 
& \leq &
\nu^1([0,1]) 
\ = \
\int_0^1 \frac{d\nu_1}{d\lambda} \, d\lambda 
\ = 
\int_{R_1} \frac{d\nu_1}{d\lambda} \, d\lambda 
\, + \int_{R_1^c} \frac{d\nu_1}{d\lambda} \, d\lambda \nonumber \\
\label{Rlbd}
& \leq & 
\int_{R_1} 1 d\lambda + \int_{R_1^c} \delta \, d\lambda
\ \leq \ 
\lambda(R_1) + \delta. 
\end{eqnarray}
As $\delta = \eta/12$, we have $\lambda(R_1) \geq \eta/12 > 0$.  This
completes the first stage of the construction.

\vskip.1in

\noindent
{\bf Further Stages.}
Subsequent stages follow the
general iterative procedure used to construct $R_1$. 
Let $\omega_{n,s}$, $f_{n,s}$, $J_{n,s}$, $m_{n,s}$, $H_{n,s}$ and $G_{n,s}$ denote
the various quantities appearing at the $n$th step of stage $s$.
In particular, let $f_{n,1} = f_n$ be the $n$'th function produced
at stage 1, and define $J_{n,1}$, $m_{n,1}$, $H_{n,1}$ and $G_{n,1}$
in a similar fashion.

Suppose that for some $s \geq 2$ the construction of the splitting
sets $R_1,\ldots,R_{s-1}$ is complete, and that we wish
to construct the set $R_s$ at stage $s$.  Let $f_{1,s}$ be any element of $\F$,
and suppose that $f_{1,s}, \ldots, f_{n,s}$ have already been
selected.  Define the join
\[ 
J_{n,s} \ = \ \D_n \vee 
\bigvee_{i=1}^n \pi(f_{i,s}) \vee \bigvee_{j=1}^{s-1} \{ R_j, R_j^c \} .
\] 
It follows from the ergodic theorem and Proposition \ref{Steele}
that there exists a sample point 
$\omega_{n+1,s} \in \Omega$, an integer 
$m_{n+1,s} \geq 1$, and a function $f_{n+1,s} \in \F$ such that 
\be
\label{inq1}
\Delta^{\omega_{n+1,s}} (I_A : m_{n+1,s}) \leq  \delta \, \lambda(A)
\mbox{ for each } A \in J_{n,s} 
\ee
and 
\be
\label{inq2}
\Delta^{\omega_{n+1,s}} (f_{n+1,s} : m_{n+1,s}) > \eta .
\ee
We may then define the join $J_{n+1,s}$ using $f_{n+1,s}$ and
continue in the same fashion.
For each $n \geq 1$ define the family
\[
H_{n,s}
\ = \ 
\left\{ A \in J_{n,s} \ : \  
\Delta^{\omega_{n+1,s}} ( f_{n+1,s} \cdot I_A : m_{n+1,s}) 
\, > \, \frac{\eta}{2} \lambda(A) \right\} 
\]
and $G_{n,s} = \cup \, H_{n,s} \subseteq [0,1]$.  Lemma \ref{Glem1} ensures that 
$\lambda(G_{n,s}) \geq \eta / 6$. 

As in stage 1, there is a sequence of integers
$n_s(1) < n_s(2) < \cdots$ such that the sub-probability measures
$\lambda_{r,s}(B) =  \lambda(B \cap G_{n_s(r),s} )$ 
converge weakly as $r \to \infty$
to a sub-probability measure $\nu_s$ on $([0,1], \sigB)$ that is absolutely
continuous with respect to $\lambda(\cdot)$.  
Define $R_s = \{ x: (d\nu_s / d\lambda)(x) > \delta \}$.  The argument
in (\ref{Rlbd}) shows that $\lambda(R_s) \geq \eta/12$.  In what follows,
we need to consider density points of $R_s$.  To this end,
for each $s \geq 1$ let 
\vskip.05in
\[
\tilde{R}_s
\ = \ 
\left\{ x \in R_s : \lim_{\alpha \to 0} 
\frac{ \lambda((x-\alpha ,x+\alpha) \cap R_s)}
                                                           {2 \alpha} \, = \, 1 \right\} .
\]
\vskip.06in
\noindent
be the Lebesgue points of $R_s$.
By standard results on differentiation of integrals (c.f.\ Theorem 31.3 of
Billingsley (1995)), we have
$\lambda (\tilde{R}_s) = \lambda(R_s) \geq \eta / 12$.

\subsection{Existence of the Intersection Tree}

Fix an integer $L \geq 1$.  As the measures of the sets $\tilde{R}_s$ 
are bounded away from 
zero, there exist positive integers $s_0 < s_1 < \ldots < s_{L}$ such that                                                                                                                                                               
$\lambda(\bigcap_{j=0}^{L} \tilde{R}_{s_j}) > 0$.    Define the intersections 
\[
Q_l = \bigcap_{j=0}^{L-l} \tilde{R}_{s_j}
\] 
for $l = 0, 1, \ldots, L$, and note that $Q_l \subseteq Q_{l+1}$.
In what follows, $B^o$, $\overline{B}$ and $\partial B$ denote,
respectively, the interior, closure and boundary of a set $B \subseteq [0,1]$.
The following result is a strengthened version of Proposition \ref{Tree1}
that incorporates the sets $Q_l$.  
Its proof completes the proof of Proposition \ref{Tree1}.

\begin{prop}
\label{Tree2}
Suppose that $\Gamma(\F : \bX) > \eta > 0$ and that 
every element of $\C(\F)$ is a finite union of intervals.  
Then there exists functions $g_1, \ldots, g_L \in \F$ and 
a complete binary tree $T$ of depth $L$ 
such that each node $t \in T$ is associated with a subset $B_t$ of $[0,1]$
subject to the following conditions:
\begin{enumerate}

\item[(a)]
For each internal node $t \in T[l]$, the
sets $B_{t'}$ and $B_{t''}$ associated with its children
$t'$ and $t''$ are 
equal to non-adjacent $\eta/5$-segments of $g_{l+1}$.

\item[(b)]
For each node $t \in T$, the intersection $W_t$ of the sets $B_s$ 
appearing along a descending path from the root to $t$ 
has non-empty interior.

\item[(c)]
If $t \in T[l]$ then the intersection $W_t^o \cap Q_l$ is non-empty.

\end{enumerate}
\end{prop}

\vskip.1in

\noindent
{\bf Proof of Proposition \ref{Tree2}}:  
Let $T$ be a complete binary tree of depth $L$ with root $t_0$, and
let $B_{t_0} = [0,1]$.  We will assign sets $B_t$ to the nodes of $T$ on 
a level-by-level basis, beginning with the children of the root.
We show below that there exists a 
function $g_1 \in \F$, and non-adjacent $\gamma$-segments $U, V \in \pi(g_1)$,
such that $U^o \cap Q_1$ and $V^o \cap Q_1$ are non-empty.
The children of $t_0$ may then be associated with $U$ and $V$, in either
order.
To begin, choose a point $x_1 \in Q_0$, which is non-empty by construction, 
and let $\epsilon = \delta / 2 (\delta + 1)$.  It follows from the 
definition of the sets $\tilde{R}_{s_j}$, that there exists  
$\alpha_1 > 0$ such that 
$I_1 \deq (x_1 - \alpha_1, x_1 + \alpha_1)$ satisfies
\be
\label{Int1prop}
\lambda(I_1 \cap Q_0) 
\ \geq \ 
(1 - \epsilon) \lambda(I_1) 
\ = \ 
2 \alpha_1 (1 - \epsilon) .
\ee
To simplify notation, let $\kappa = s_L$.  
The last display and the definition of $R_\kappa$ imply that
\[
\nu_{\kappa}(I_1 \cap R_{\kappa})
\ = \ 
\int_{I_1 \cap R_\kappa} \frac{ d\nu_{\kappa} }{d\lambda} \, d\lambda
\ > \ 
\delta \, \lambda(I_1 \cap R_\kappa)
\ \geq \ 
2 \alpha_1 (1 - \epsilon) \delta .
\]
Let $\{ n_\kappa(r) : r \geq 1 \}$ be the subsequence used to define
the sub-probability $\nu_\kappa$. 
As $I_1$ is an open set, it follows from the Portmanteau theorem that
\[
\liminf_{r \to \infty} 
\lambda( I_1 \cap G_{n_\kappa(r),\kappa} )
\ \geq \ 
\nu_{\kappa} (I_1)
\ \geq \ 
\nu_{\kappa} (I_1 \cap R_\kappa)
\ > \ 
2 \alpha_1 (1 - \epsilon) \delta .
\]
Choose $r$ sufficiently large so that 
$\lambda( I_1 \cap G_{n_\kappa(r),\kappa} ) ) \, > \, 
2 \alpha_1 (1 - \epsilon) \delta$ and
$2^{ - n_\kappa(r) } < \delta \, \alpha_1 / 4$.  We require
the following subsidiary lemma.  Its proof is identical to
Lemma 4 in \cite{AdNob09}, but is included in the Appendix for
completeness.

\vskip.16in
 
\begin{lem}
\label{Alem}
There exists a set $A \in H_{n_\kappa(r),\kappa}$ such that
$A \subseteq I_1$ and $\lambda(A \cap Q_1) > 0$.  Moreover,
$A$ is contained in $Q_1$. 
\end{lem}

\vskip.1in

Let $g_1 = f_{ n_{\kappa}(r)+1, \kappa } \in \F$.  By assumption, 
each element of $\pi(g_1)$ is a 
finite union of intervals, and no random variable $X_i$ takes values in the 
finite set $\cup_{C \in \pi(g_1)} \partial C$. 
We argue that the set $A$ identified in Lemma \ref{Alem}
(and therefore $Q_1$) has non-empty intersection
with the interiors of two non-adjacent segments of $g_1$.  
As $A$ has positive measure, and the boundary of each segment 
of $g_1$ has measure zero, it suffices to exclude the possibility that $A$ 
intersects no segments, 
only one segment, or only two adjacent segments of $g_1$.  

As $\lambda(A) > 0$ and the segments of $g_1$ form a partition of $[0,1]$,
$A$ must intersect
the interior of at least one segment of $g_1$.   
Suppose that $A$ intersects only one segment $U = s_{k}(g_1)$ of $g_1$.
Let $h(x) = g_1(x) - (k - 1) \gamma$, and note that $0 \leq h(x) \leq \gamma$
for each $x \in U$.    In this case,  
\begin{eqnarray}
E(g_1 \, I_A)(X) 
& = & 
\sum_{C \in \pi(g_1)} E(g_1 \, I_A \, I_C)(X) 
\ = \ 
E(g_1 \, I_A \, I_U)(X) \nonumber \\ [.1in]
\label{bndexp}
& = &
\gamma (k -1) \lambda(A) \, + \, E(h \, I_A)(X) .
\end{eqnarray}
Similarly, for each $m \geq 1$, 
\begin{eqnarray}
\frac{1}{m} \sum_{i=1}^m (g_1 \, I_A)(X_i) 
& = & 
\frac{1}{m} \sum_{i=1}^m \sum_{C \in \pi(g_1)} (g_1 \, I_A \, I_C)(X_i) 
\ = \ 
\frac{1}{m} \sum_{i=1}^m (g_1 \, I_A \, I_U)(X_i) \nonumber \\ [.1in]
\label{bndsum}
& = & 
\gamma (k -1) \, \frac{1}{m} \sum_{i=1}^m I_A (X_i) \ + \  
\frac{1}{m} \sum_{i=1}^m (h \, I_A)(X_i) .
\end{eqnarray}
Letting $m = m_{n_\kappa(r) + 1, \kappa}$, we find that
\begin{eqnarray*}
\frac{\eta}{2} \, \lambda(A) 
& < &
\Delta^w(g_1 \cdot I_A  : m) \\
& \leq &
\gamma (k -1) \, \Delta^w(I_A  : m)
\, + \,
\max\left\{ \frac{1}{m} \sum_{i=1}^m (h \, I_A)(X_i), \, E(h \, I_A)(X) \right\} \\ [.1in]
& \leq &
\gamma (k -1) \, \Delta^w(I_A  : m)
\, + \,
\gamma \max\left\{ \frac{1}{m} \sum_{i=1}^m I_A (X_i), \, \lambda(A) \right\} \\ [.1in]
& \leq &
\gamma (k -1) \, \Delta^w(I_A  : m)
\, + \,
\gamma (\lambda(A) + \Delta^w(I_A  : m)) \\ [.05in]
& \leq &
\Delta^w(I_A  : m) + \gamma \lambda(A) \\ [.05in]
& \leq &
(\delta + \gamma) \, \lambda(A) .
\end{eqnarray*}
Here the first inequality follows from the definition of 
$H_{n_\kappa(r), \kappa}$, the second follows 
from (\ref{bndexp}) and
(\ref{bndsum}), the third follows from the bound
on $h(\cdot)$, and last follows 
from the definition of $m$.  Comparing the first and last
terms above, our definition of $\delta = \eta / 12$ and 
$\gamma = \eta / 5$ yields a contradiction.

Suppose finally that $A$ intersects only two adjacent segments of $g_1$, 
say $U = s_{k}(g_1)$ and $V = s_{k+1}(g_1)$.
Let $h(x)$ be defined as above, and 
note that $0 \leq h(x) \leq 2 \gamma$ for $x \in U \cup V$.  
Arguing as above, we find that
\[
E(g_1 \cdot I_A)(X) 
\ = \
\gamma (k -1) \lambda(A) \, + \, E(h \, I_A)(X) ,
\]
and that for each $m \geq 1$, 
\[
\frac{1}{m} \sum_{i=1}^m (g_1 \, I_A)(X_i) 
\ = \ 
\gamma (k -1) \, \frac{1}{m} \sum_{i=1}^m I_A (X_i) \ + \  
\frac{1}{m} \sum_{i=1}^m (h \, I_A)(X_i) .
\]
Letting $m = m_{n_\kappa(r) + 1, \kappa}$, the previous
two displays, and arguments like those above, can be 
used to show that
\begin{eqnarray*}
\frac{\eta}{2} \, \lambda(A) 
& < &
\Delta^w(g_1 \cdot I_A  : m) \\
& \leq &
\gamma (k -1) \, \Delta^w(I_A  : m)
\, + \,
2 \gamma (\lambda(A) + \Delta^w(I_A  : m)) \\
& \leq &
(1 + \gamma) \, \Delta^w(I_A  : m) + 2 \gamma \lambda(A) \\
& \leq &
((1 + \gamma) \delta + 2 \gamma) \, \lambda(A) .
\end{eqnarray*}
Comparing the first and last terms, the definition of 
$\delta = \eta / 12$ and $\gamma = \eta/5 $ yields a contradiction,
and we conclude that $A$ intersects the interiors of two non-adjacent
segments $U$ and $V$ of $g_1$.  This completes the assignment of
sets to the children of the root $t_0$.

\vskip.1in

Suppose now that for some $l \leq L-1 $ we have assigned sets
$B_t \subseteq [0,1]$ to each node $t$ of $T$ having depth less
than or equal to $l$, in such a way that properties (a) - (c) of
the Proposition hold.  There are $2^l$ nodes of $T$ at distance $l$ from the
root.  Denote these nodes by $1 \leq j \leq 2^l$, and let $W_j$ be the intersection
of the sets $B_s$ appearing on the descending path from the root $t_0$ of
$T$ to node $j$ at level $l$.  By assumption,
$W_j^o \cap Q_l$ is non-empty:
let $x_j \in W_j^o \cap Q_l$ for each $j \in [2^l]$.
Select $\alpha_{l+1} > 0$ such that, for each $j$, the interval
$I_j \deq (x_j - \alpha_{l+1}, x_j + \alpha_{l+1})$
is contained in $W_j^o$ and satisfies
\[
\lambda(I_j \cap Q_l) 
\ \geq \ 
(1 - \epsilon) \lambda(I_j) 
\ = \
2 \alpha_{l+1} (1 - \epsilon) . 
\]
Let $\kappa' = s_{L - l}$ and let $\{ n_{\kappa'} (r) : r \geq 1 \}$ 
be the subsequence used to define the sub-probability $\nu_{\kappa'}$. 
For each interval $I_j$, 
\[
\liminf_{r \to \infty} 
\lambda( I_j \cap G_{n_{\kappa'} (r), \kappa'} )
\ \geq \ 
\nu_{\kappa'} (I_j)
\ \geq \ 
\nu_{\kappa'} (I_j \cap R_{\kappa'})
\ > \ 
2 \alpha_{l+1} (1 - \epsilon) \delta .
\]
where the last inequality follows from the previous display,
and the fact that $Q_l \subseteq R_{\kappa'}$. 
Choose $r$ sufficiently large so that 
$\lambda( I_j \cap G_{n_{\kappa'} (r), \kappa'} ) >
2 \alpha_{l+1} (1 - \epsilon) \delta$ for each $j = 1,\ldots,2^l$,
and $2^{ - n_{\kappa'}(r) } < \delta \, \alpha_{l+1} / 4$.  

Applying the proof of Lemma \ref{Alem} to each interval $I_j$,
we may identify sets 
$A_1, A_2, \ldots, A_{2^l} \in H_{n_{\kappa'} (r), \kappa'}$ such that 
$\lambda(A_j) > 0$,
$A_j \subseteq I_j \subseteq W_j^o$, and
$A_j \subseteq  Q_{l+1}$ for each $j = 1,\ldots,2^l$.
Define $g_{l+1} = f_{n_{\kappa'} (r)+1, \kappa'} \in \F$.  
Arguments
identical to those in the case $l=0$ above show that, for each $j$, there
exist non-adjacent segments $U_j,V_j$ of $g_{l+1}$ such that
$A_j \cap U_j^o$ and $A_j \cap V_j^o$ are non-empty.   Assigning
the sets $U_j$ and $V_j$ to the left and right children of $j$ in $T$, in either order,
ensures that property (a) of the proposition is satisfied.  
For the child $t$ of node
$j$ associated with the set $U_j$ we have $W_t = W_j \cap U_j$.  It follows from the fact that
$A_j \subseteq W_j^o$, $A_j \cap U_j^o \neq \emptyset$ and 
$A_j \subseteq  Q_{l+1}$ that $W_t^o \cap Q_{l+1} \neq \emptyset$, and therefore
properties (b) and (c) of the proposition are satisfied.  
The argument for the other child of node $j$ is similar.  
This completes the proof of Proposition \ref{Tree2}.

\section{Proof of Proposition \ref{FullJoin}}
\label{PfFullJoin}

\noindent
{\bf Proof of Proposition \ref{FullJoin}:}  
Fix $L \geq 1$ such that $2^{L-1} / K^2 \geq 4$, and 
let $T$ be the complete binary tree of depth $L$ described in Proposition
\ref{Tree1}.  Suppose that each interior node in $t \in T$ is labeled with 
the indices of the segments assigned to its children: if
the segments $s_k(g_r)$ and $s_{k'}(g_r)$ of $g_r$ are assigned
to the children of a node $t \in T[r-1]$, then 
$t$ is assigned the label $\ell(t) = (k,k') \in [K]^2$, where
$[K] = \{1,\ldots, K\}$.

Let $L_0 = L-1$.  
By an elementary pigeon-hole argument, there exist non-adjacent
integers $k_0,k_0' \in [K]$ such that the set $S_0$ of nodes
$t \in T[L_0]$ with $\ell(t) = (k_0,k_0')$ has cardinality at least
$2^{L_0} / K^2$.  (Here $K^2$ is an upper bound on the number 
of non-adjacent pairs $k, k' \in [K]$.)  
Let $u_0 = \lceil \log_2 K^2 + 1 \rceil$. 

It follows from Lemma \ref{Ptree} and an additional pigeon hole 
argument that there exists an integer $L_1$, a pair
$k_1,k_1' \in [K]$ of non-adjacent integers,
and a set of nodes $S_1 \subseteq T[L_1]$ with the following properties:
(i) $L_0 - u_0 \leq L_1 \leq L_0 - 1$; 
(ii) $\ell(t) = (k_1,k_1')$ for every $t \in S_1$;
(iii) for every $t \in S_1$, each child of $t$ is an ancestor of $S_0$; and
(iv) $|S_1|  \ \geq \ 2^{L_0} / 4 L K^4$.
In particular, inequalities (i) and (iv) imply that
\be
\label{scardineq}
|S_1| 
\ \geq \ 
2^{L_1} \left( \frac{ 2^{L_0 - L_1} }{ 4 L K^4 } \right)
\ \geq \ 
2^{L_1} \left( \frac{ 1 }{ 2 L K^4 } \right)
\ \geq \ 
\frac{ 2^{L_0} }{ 8 L K^6 } .
\ee
If the last term above is greater than or equal to $4$, then we may
apply Lemma \ref{Ptree} again to find an integer $L_2$ and a set of nodes 
$S_2 \subseteq T[L_2]$ with properties
analogous to (i) - (iv) above.
Continuing in this fashion, we obtain integers $L_0 > L_1 > \cdots > L_R \geq 0$, 
sets of nodes $S_r \subseteq T[L_r]$, and non-adjacent pairs $k_r,k_r' \in [K]$ 
such that for $1 \leq r \leq R$ and for
every node $t \in S_r$, $\ell(t) = (k_r,k_r')$ and both children 
of $t$ are ancestors of $S_{r-1}$.
In particular, using arguments like those in (\ref{scardineq}),
one may show that 
\[
|S_r| 
\ \geq \ 
2^{L_r} 
\left( \frac{ 1 }{ (2 L K^2)^r K^2 } \right)
\ \geq \ 
\frac{ 2^{L - 1} }{ 4^r \cdot K^{2r + 1} \cdot (2 L K^2)^{r(r+1) / 2} } ,
\]
and therefore $R = R(L)$ can be taken to be the largest integer $r \geq 1$
for which the last term above is greater than $4$. 
In particular, $R(L)$ tends to infinity with $L$.

From the construction above, and an additional pigeon-hole argument,
we may identify an integer 
$N = N(L) \geq R(L)/K^2$ and a subsequence $i_0 < i_1 < \cdots < i_N$ 
of $L_R, L_{R-1},\ldots, L_0$ such that
$(k_{i_j},k_{i_j}') = (k,k')$ for a fixed non-adjacent pair $(k,k') \in [K]^2$.
From the associated node-sets $S_{i_0},\ldots, S_{i_N}$ one may construct an
embedded binary subtree $T_o$ of $T$ all of whose node labels are equal
to $(k,k')$.  To see this,
let the root of $T_o$ be any node $s \in S_{i_0}$.  At each level
$0 \leq r \leq N-1$ let
the left and right children of $t \in T_o[r]$ be (necessarily distinct) descendants 
in $S_{i_r + 1}$ of the children of $t \in T$.  Then it is easy
to see that $T_o$ is a complete binary
tree of depth $N$.

For $r = 0,\ldots,N-1$ let $h_r = g_{i_r + 1}$.
By construction, each node $t \in T_o[r]$ is contained in
$S_{i_r}$ and has label $\ell(t) = (k,k')$.  Thus 
the children $t'$ and $t''$ of $t$ in $T_o$ are associated with the segments
$s_k(h_r)$ and $s_{k'}(h_r)$ of $h_r$.  For each terminal node $t \in \tilde{T}_o$ 
let $W_t$ be the intersection 
of the sets $B_s$ appearing on the descending path (in $T$) from
the root of $T_o$ to $t$.  The construction of $T_o$ ensures that every member
of $\{ W_t : t \in \tilde{T}_o \}$ is contained in a unique element of the join
\[
J \ = \ \{ s_k(h_0), s_{k'}(h_0) \} \vee \cdots \vee \{ s_l(h_{N-1}), s_{l'}(h_{N-1}) \}
\]
Moreover, by Proposition \ref{Tree1}, each set $W_t$ has non-empty
interior, and positive Lebesgue measure, and the same is therefore true 
for each element of $J$.
As $N(L)$ tends to infinity with $L$, the lemma follows.

\section{Proof of Theorem \ref{FS2}} 
\label{PfFS2}

\noindent
{\bf Proof of Theorem \ref{FS2}:} 
Let $\F$ and $\bX$ be as in the statement of the 
proposition.  Then
$\Gamma(\F : \bX) > \eta > 0$.
Let $\C(\F)$ be the countable family defined in (\ref{regsets}).
Without loss of generality, we may assume
that $\F$ contains the identity function $f_0(x) = x$,
and therefore $\C(\F)$ satisfies the shrinking
diameter condition of Lemma \ref{RedLem3}.
Let the sets $V_1,V_2 \subseteq [0,1]$ and map 
$\phi(\cdot)$ be as in the
statement of Lemma \ref{RedLem3}.  

Define random variables $Y_i = \phi(X_i)$ for $i \geq 1$.  Then the process
$\bY = Y_1, Y_2, \ldots$ is stationary and ergodic with $Y_i \sim \lambda$.
For each $f \in \F$ define an associated function $g_f: [0,1] \to [0,1]$ via
the rule
\be
\label{gfdefn}
g_f (u) \ = \ 
\left\{ \begin{array}{ll} 
          (f \circ \phi^{-1}) (u) & \mbox{if $u \in V_2$} \\ 
          0 & \mbox{if $u \in V_2^c$} 
               \end{array} \right. 
\ee
and let ${\cal G} = \{ g_f : f \in \F \}$.  Arguments like those in Section \ref{RUI}
above show that $\Gamma_m({\cal G} : \bY )$ is equal to $\Gamma_m(\F : \bX)$ with
probability one for each $m \geq 1$, and 
consequently $\Gamma({\cal G} : \bY) > \eta$.

Let the constants $\gamma$ (equal to $\eta / 5$) and $K$, and the segments $s_k(f)$, 
be defined as in (\ref{segdef}), and let
$\epsilon = \Gamma({\cal G} : \bY) - \eta > 0$.
Choose a finite sequence of rational numbers 
$0 = a_0 < a_1 < \cdots < a_N = 1$ that includes
$\{ \gamma k : k= 1, \ldots, K-1 \}$ and is such that
$\max_j |a_j - a_{j-1}| < \epsilon / 2$.
Define intervals $U_j = [a_{j-1},a_j)$
for $j = 1,\ldots, N-1$, and let $U_N = [a_{N-1},1]$.  Using (\ref{gfdefn})
one may verify  
that for each $g_f \in {\cal G}$,
\[
g_f^{-1} U_j
\ = \ 
\left\{ \begin{array}{ll} 
          \phi(f^{-1} U_j) & \mbox{if $2 \leq j \leq N$} \\ 
          \phi(f^{-1} U_j) \cup V_2^c & \mbox{if $j = 1$} 
               \end{array} \right.,
\]
where the second condition results from the fact that the interval $U_1$ contains zero.

Let ${\cal U}$ be the family of subsets of $[0,1]$ that are equal to
finite unions of intervals, and let $A \cong B$ denote the fact that $A$ and $B$ are
equivalent mod 0, in other words, $\lambda(A \triangle B) = 0$.  Fix a function $f \in \F$,
and let $g_f$ be the associated element of ${\cal G}$.
Lemma \ref{RedLem3} and the fact that $\lambda(V_2^c) = 0$ 
imply that there exists sets $C_1, \ldots, C_N \in {\cal U}$ 
such that $g_f^{-1} U_j \cong C_j$ for $1 \leq j \leq N$.  If $i \neq j$ then 
\[
\lambda(C_i \cap C_j) 
\ = \ 
\lambda(g_f^{-1} U_i \cap g_f^{-1} U_j)
\ = \ 
\lambda(g_f^{-1} (U_i \cap U_j))
\ = \ 
0 
\]  
so that $C_i$ and $C_j$ can intersect only at 
the endpoints of their constitutive intervals.  
It follows that the function $h_f(u) = \sum_{j=1}^N a_{j-1} I_{C_j}(u)$
approximates $g_f$ in the sense that $|g_f(u) - h_f(u)| < \epsilon / 2$ 
with probability one.  Moreover, $h_f^{-1} [a,b) \in {\cal U}$ 
for all rational $a, b$.  

Let ${\cal H} = \{ h_f : f \in \F \}$ be the family of simple approximations to 
the elements of $\G$.  Then $\C({\cal H})$ is contained in ${\cal U}$, and
a straightforward argument shows that $\Gamma({\cal H} : \bY) > \eta$.
Fix $L \geq 1$.  
As ${\cal H}$ satisfies the conditions of
Proposition \ref{FullJoin}, there exist functions $f_1,\ldots,f_L \in \F$ and a pair
of non-adjacent integers $k,k' \in [K]$ such that the join
\[
J_h = \bigvee_{\ell=1}^L \{ s_k(h_{f_\ell}), s_{k'}(h_{f_\ell}) \}
\]
has $2^L$ elements, each with positive measure.
In order to obtain a full join for the segments of $f_1,\ldots,f_L$,  
we examine how the segments of $h_f$ are related to those of $f$.
To this end, let $i < j$ be such that
$a_i =(k-1)\gamma$ and $a_j = k \gamma$.  
Then for every $f \in \F$,
\begin{eqnarray*}
s_k(h_f)
& = &
h_f^{-1} [(k-1) \gamma, k \gamma) \
\ = \ 
h_f^{-1} [a_i, a_j) \\ [.05in]
& = & 
\bigcup_{r=i}^{j-1} C_{r+1} 
\ \cong \ 
\bigcup_{r=i}^{j-1} g_f^{-1} U_{r+1} \\ [.03in]
& = &
g_f^{-1} [(k-1) \gamma, k \gamma) \\ [.04in]
& \cong &
\phi( f^{-1} [(k-1)\gamma, k \gamma) ) 
\ = \ \phi( s_k(f) ) .
\end{eqnarray*}
The same argument applies to $s_{k'}(h_f)$, and therefore
every element of $J_h$ is equivalent mod zero to an element
of the join
\[
J_h'
\ = \
\bigvee_{\ell=1}^L \{ \phi(s_k(f_\ell)) , \phi(s_{k'}(f_\ell)) \} .
\]
As $\phi$ is a bijection almost everywhere,
every element of $J_h'$ is equivalent mod zero
to a set of the form $\phi(A)$, where $A$ is an element of the join
\[
J_f \ = \ \bigvee_{\ell=1}^L \{ s_k(f_\ell) , s_{k'}(f_\ell) \} .
\]
As each cell of $J_h'$ has positive Lebesgue measure, the same
is true of the cells of $J_f$.  In particular, 
$J_f$ has (maximum) cardinality $2^L$.  As $L \geq 1$
was arbitrary, Theorem \ref{FS2} follows from 
Lemma \ref{Join}.

\appendix

\section{Appendix}
\label{Tech}

\subsection{Proof of Lemma \ref{Alem}}

The proof of Lemma \ref{Alem} appears in \cite{AdNob09}; we reproduce it
here for completeness.

\noindent
{\bf  Proof:} Let $G = G_{n_\kappa(r),\kappa}$.  The  
choice of $n_\kappa(r)$ ensures that 
\begin{eqnarray*}
(1 - \epsilon) \, \delta \, \lambda(I_1)
& \leq &
\lambda( I_1 \cap G ) \\
& = & 
\lambda( I_1 \cap Q_1 \cap G )
\ + \ 
\lambda( I_1 \cap Q_1^c \cap G ) \\
& \leq &
\lambda( I_1 \cap Q_1 \cap G )
\ + \ 
\lambda( I_1 \cap Q_1^c ) \\
& \leq &
\lambda( I_1 \cap Q_1 \cap G )
\ + \ 
\epsilon \lambda(I_1) 
\end{eqnarray*}
where the final inequality follows from (\ref{Int1prop}) and the fact
that $Q_0 \subseteq Q_1$.
It follows from the display and the definition of $\epsilon$
that $\lambda( I_1 \cap Q_1 \cap G ) \geq \delta \alpha_1$.
As the collection of sets used to define the join $J_{n_\kappa(r),\kappa}$ 
includes the dyadic intervals of order $n_\kappa(r)$, each element $A$
of the join has diameter (and Lebesgue measure) bounded by 
$2^{ - n_\kappa(r) } < \delta \, \alpha_1 / 4$.  These last two inequalities
imply that
\[
\delta \, \alpha_1
\ \leq \ 
\lambda( I_1 \cap Q_1 \cap G )
\ \leq \ 
\sum_{A} \lambda(Q_1 \cap A) \, + \,  2 \, \frac{\delta \, \alpha_1}{4} ,
\]
where the sum is over $A \in H_{n_\kappa(r),\kappa}$ such that
$A \subseteq I_1$.  In particular, it is clear that the sum is necessarily 
positive, and the first part of the claim follows.  Moreover, for
any set $A \in H_{n_\kappa(r),\kappa}$ the definition
of the join $J_{n_\kappa(r),\kappa}$ requires that $A$ be contained 
in either $R_{k_j}$ or $R_{k_j}^c$ for each $j=0,\ldots,L-1$.
If $\lambda(A \cap Q_1) > 0$ then necessarily
$A \cap Q_1 \neq \emptyset$, and these containment relations
imply that $A \subseteq Q_1$.  This completes the proof of
Lemma \ref{Alem}

\vskip.3in

\noindent
{\bf\Large  Acknowledgements} \\
The work presented in this paper was supported in part
by NSF grant DMS-0907177.

\end{document}